\documentclass[a4paper,11pt]{article}
\usepackage[margin=1in]{geometry}
\usepackage{tikz}
\usepackage[small,bf,hang]{caption}
\usepackage{float}
\usepackage{siunitx}
\usepackage{amsmath,amsthm,amssymb}
\usepackage{tikz}
\usepackage{hyperref}
\hypersetup{colorlinks=true, citecolor = blue, linkcolor=black, filecolor=magenta, urlcolor=cyan}
  
\newtheorem{theorem}{Theorem}

\newtheorem{corollary}{Corollary}
\newtheorem{proposition}{Proposition}
\newtheorem{lemma}{Lemma}

\newtheoremstyle{case}{}{}{}{}{}{:}{ }{}
\theoremstyle{case}

\title{Recoloring some hereditary graph classes}

\author{Manoj Belavadi
 \thanks{Department of Mathematics, Wilfrid Laurier University,
 Waterloo, ON, Canada, N2L 3C5. Email: \texttt{mbelavadi@wlu.ca}. ORCID: 0000-0002-3153-2339. Research supported by the Natural Sciences and Engineering Research Council of Canada (NSERC) grant RGPIN-2016-06517.}
 \and Kathie Cameron
 \thanks{Department of Mathematics, Wilfrid Laurier University,
 Waterloo, ON, Canada, N2L 3C5. Email: \texttt{kcameron@wlu.ca}. ORCID: 0000-0002-0112-2494. Research supported by the Natural Sciences and Engineering Research Council of Canada (NSERC) grant RGPIN-2016-06517.}}

\begin{document}
\maketitle

\begin{abstract}
The reconfiguration graph of the $k$-colorings, denoted $R_k(G)$, is the graph whose vertices are the $k$-colorings of $G$ and two colorings are adjacent in $R_k(G)$ if they differ in color on exactly one vertex. A graph $G$ is said to be recolorable if $R_{\ell}(G)$ is connected for all $\ell\geq \chi(G)$+1. In this paper, we study the recolorability of several graph classes restricted by forbidden induced subgraphs. We prove some properties of a vertex-minimal graph $G$ which is not recolorable. We show that every (triangle, $H$)-free graph is recolorable if and only if every (paw, $H$)-free graph is recolorable. Every graph in the class of $(2K_2,\ H)$-free graphs, where $H$ is a 4-vertex graph except $P_4$ or $P_3$+$P_1$, is recolorable if $H$ is either a triangle, paw, claw, or diamond. Furthermore, we prove that every ($P_5$, $C_5$, house, co-banner)-free graph is recolorable.
\end{abstract}

\section{Introduction}
Reconfiguration problems are defined on various source problems in the literature. Typical questions asked are: given any two solutions to a source problem can we transform one to another by a sequence of specified elementary steps and find the length of such a sequence if it exists. Reconfiguration problems have applications in many fields such as combinatorial games, Glauber dynamics, and cellular networks. Reconfiguration problems have been defined on many problems in graph theory such as vertex coloring, independent set, dominating set, and matching. See a survey by Nishimura \cite{Nishimura} for more information on various reconfiguration problems. One of the ways to study reconfiguration problems is to define the reconfiguration graph. The reconfiguration graph of a source problem is the graph whose vertices corresponds to the solutions of the source problem and an edge corresponds to an elementary step that transforms one solution to another. We can then ask: Given any two vertices in the reconfiguration graph, is there a path between those two vertices? What is the length of such a path? Is the reconfiguration graph connected? What is the diameter of the reconfiguration graph? In this paper we study the reconfiguration graph of vertex colorings, its connectivity and its diameter.

Let $G$ be a finite simple graph with vertex-set $V(G)$ and an edge-set $E(G)$. We use $n$ to represent the number of vertices in a graph. Two vertices $u$ and $v$ are \emph{adjacent} in $G$ if $uv\in E(G)$. A \emph{path} is a sequence of distinct vertices $v_i$ and edge $v_iv_j$ of the form $v_o$, $v_ov_1$, $v_1,\dots, v_{p-1}$, $v_{p-1}v_p$, $v_p$. The length of a path is equal to the number of edges in the path. A graph is said to be \emph{connected} if there exists a path between every pair of distinct vertices of the graph.  The \emph{diameter} of a graph is the length of a longest shortest path between any two distinct vertices of the graph. The  subgraph of $G$ \emph{induced} by a subset $S\subseteq V(G)$ is the graph whose vertex-set is $S$ and whose edge-set is all edges of $G$ which join vertices in $S$; it is denoted by $G[S]$. For a positive integer $k$, a $k$-coloring of $G$ is a mapping from $V(G)$ to a set of colors $\{1,2,\dots,k\}$ such that no pair of adjacent vertices receive the same color. We say that $G$ is $k$-colorable if it admits a $k$-coloring, and the \emph{chromatic number} of $G$, denoted $\chi(G)$, is the minimum number of colors required to color $G$. The \emph{reconfiguration graph of the $k$-colorings}, denoted $R_k(G)$, is the graph whose vertices are the $k$-colorings of $G$ and two colorings are adjacent in $R_k(G)$ if they differ in color on exactly one vertex. A graph $G$ is said to be \emph{recolorable} if $R_{\ell}(G)$ is connected for all $\ell\geq \chi(G)$+1. The $\ell$ \emph{recoloring diameter} of a graph $G$ is the diameter of the reconfiguration graph $R_{\ell}(G)$ where $\ell\geq \chi(G)$+1. For a graph $H$, a graph $G$ is $H$-free if no induced subgraph of $G$ is isomorphic to $H$. For a collection of graphs $\mathcal{H}$, a graph $G$ is $\mathcal{H}$-free if $G$ is $H$-free for every $H \in \mathcal{H}$. Let $P_n$, $C_n$, and $K_n$ denote the path, cycle, and complete graph on $n$ vertices, respectively. Let $K_{n,n}$ denote the complete bipartite graph where each partite set contain $n$ vertices. For two vertex-disjoint graphs, $G$ and $H$, the \emph{disjoint union} of $G$ and $H$, denoted by $G+H$, is the graph with vertex-set $V(G) \cup V(H)$ and edge-set $E(G) \cup E(H)$. For a positive integer $r$, we use $rG$ to denote the graph consisting of the disjoint union of $r$ copies of $G$.

The problem of  deciding if a graph $G$ is $k$-colorable for any $k\geq 3$ is NP-complete. This led many researchers to study coloring for restricted graph classes, including classes of $H$-free graphs. Deciding whether there exists a path between two colorings in $R_k(G)$ is PSPACE-complete for $k > 3$ \cite{Bonsma2009}, and can be decided in polynomial time for $k\leq 3$ \cite{Cereceda2008}. The problem remains PSPACE-complete for graph classes with bounded treewidth or bounded bandwidth \cite{Wrochna}. So it was natural to study the reconfiguration problem for $H$-free graphs.

A frozen $k$-coloring of $G$ is a $k$-coloring of $G$ where for any vertex $v\in V(G)$, each of the $k$ colors is mapped to $v$ or to a vertex adjacent to $v$. To prove that $R_k(G)$ is disconnected, we can exhibit a frozen $k$-coloring of $G$ (which corresponds to an isolated vertex in $R_{k}(G)$). Since the complete graph on $k$ vertices, $K_k$, admits a frozen $k$-coloring, it is common to study $R_{\ell}(G)$ for $\ell\geq \chi(G)$+1.

A recent question in this field was: For which graphs $H$ is every $H$-free graph recolorable? The problem was first solved for 3-vertex graphs, namely 3$K_1$, $P_2$+$P_1$, $P_3$, and $K_3$. In particular, it was proved that the classes of 3$K_1$-free, $P_2$+$P_1$-free, and $P_3$-free graphs are recolorable \cite{3k1free, bonamy2018, bonamy2014} and that there exist triangle-free graphs that are not recolorable \cite{Cereceda2008}. In \cite{bonamy2018}, Bonamy and Bousquet proved that every $P_4$-free graph $G$ is recolorable and the $\ell$ recoloring diameter is at most 4$n$. In \cite{2K2free}, Feghali and Merkel proved that there exist 2$K_2$-free graphs that are not recolorable. This question was completely answered in \cite{Belavadi}. A graph $G$ is said to be $\ell$-\emph{mixing} if $R_{\ell}(G)$ is connected.

\begin{theorem}[\cite{Belavadi}]
    Every $H$-free graph $G$ is $\ell$-mixing for all $\ell\geq \chi(G)$+1 if and only if $H$ is an induced subgraph of $P_4$ or $P_3$+$P_1$.
\end{theorem}

In \cite{Belavadi}, it was also proved that every (2$K_2$, $C_4$)-free graph is recolorable and the $\ell$ recoloring diameter is at most 4$n$. In this paper, we focus on recolorability for the class of ($H_1$, $H_2$)-free graphs, where $H_1$ and $H_2$ are not induced subgraphs of $P_4$ or $P_3$+$P_1$. In Section \ref{sec:traingle}, we investigate the recolorability of several subclasses of triangle-free graphs and give upper bounds for the diameter of the reconfiguration graph. There are 11 graphs on 4 vertices, including $P_4$, $P_3$+$P_1$, and 2$K_2$. The others are 4$K_1$, $C_4$, $K_4$, claw, co-claw, diamond, co-diamond, and paw. Since $C_6$ admits a frozen 3-coloring and is (4$K_1$, $C_4$, $K_4$, claw, co-claw, diamond, co-diamond, paw)-free, it follows that if every ($H_1$, $H_2$)-free graph $G$ is recolorable, then either $H_1$ or $H_2$ is isomorphic to 2$K_2$. Hence, in Section \ref{sec:2K2}, we study several subclasses of $2K_2$-free graphs. We prove the following results.

\begin{theorem}
    Every (2$K_2$, $H$)-free graph is recolorable when $H \in$ (triangle, paw, claw, diamond).
\end{theorem}

\begin{theorem}
    For all $p\geq 1$, there exists a $k$-colorable (2$K_2$, $H$)-free graph $G$ that is not ($k$+$p$)-mixing when $H \in$ ($4K_1$, co-diamond, co-claw).
\end{theorem}

The class of  $P_4$-sparse graphs is a superclass of $P_4$-free graphs that was first defined by Ho\`ang \cite{Hoang1985}: A graph $G$ is $P_4$-sparse if every set of five vertices induces at most one $P_4$. Biedl, Lubiw, and Merkel \cite{biedl2021} proved that every $P_4$-sparse graph $G$ is recolorable with $\ell$ recoloring diameter at most 4$n^2$. We generalize this result by proving that every ($P_5$, $C_5$, house, co-banner)-free graph $G$ is recolorable with $\ell$ recoloring diameter at most 2$n^2$.

\section{Preliminaries}
\label{sec:pre}

Here we define the terminology that is used in this paper. For a vertex $v \in V(G)$, the \emph{open neighborhood} of $v$, denoted $N(v)$, is the set of vertices adjacent to $v$ in $G$. The \emph{closed neighborhood} of $v$, denoted $N[v]$, is the set of vertices adjacent to $v$ in $G$ together with $v$. The \emph{degree of a vertex} $v$ in $G$, denoted $d_{G}(v)$, is the number of vertices in $N(v)$, and we use $d(v)$ when the context is clear. Symbols $\delta(G)$ and $\Delta(G)$ denote the minimum and the maximum degree of a vertex in $G$, respectively. A graph $G$ is said to be \emph{$p$-regular} if the degree of every vertex in $G$ is $p$. For $X,Y \subseteq V(G)$, we say that $X$ is \emph{complete} to $Y$ if every vertex in $X$ is adjacent to every vertex in $Y$. If no vertex of $X$ is adjacent to a vertex of $Y$, we say that $X$ is \emph{anticomplete} to $Y$. A \emph{component} of $G$ is a maximal connected subgraph of $G$. A \emph{join} of two graphs $G_1$ and $G_2$ is the graph obtained from the disjoint union of $G_1$ and $G_2$ by joining each vertex of $G_1$ to every vertex of $G_2$. For $p\geq 1$, we use [$p$] to denote the set $\{1,\dots,p\}$.

The set of vertices that are assigned the same color in a coloring is called a \emph{color class}. Two colorings $\alpha$ and $\beta$ of $G$ are \emph{isomorphic} if they induce the same partition of $V(G)$ into color classes. For a coloring $\alpha$ of $G$ and $X \subseteq V(G)$, we say that the color $c$ \emph{appears} in $X$ if $\alpha(x) = c$ for some $x \in X$, and we use $\alpha(X)$ for the set of colors appearing in $X$ and $|\alpha(X)|$ for the number of colors appearing on $X$. Two vertices which are not adjacent are called \emph{independent vertices} and a set of mutually independent vertices is called an \emph{independent set}. A \emph{clique} of a graph $G$ is a set of mutually adjacent vertices of $G$ and the \emph{clique number} is the size of the maximum clique in $G$. A \textit{clique cutset} $Q$ of a graph $G$ is a clique in $G$ such that $G$-$Q$ has more components than $G$. A clique cutset $Q$ of a graph $G$ is called a \textit{tight clique cutset} if there exists a component $H$ of $G$-$Q$ which is complete to $Q$. If $G$ is disconnected, then it is obvious that $G$ is recolorable if and only if every component of $G$ is recolorable, so we may assume that $G$ is connected when appropriate. Two edges are called \emph{independent} if they do not share an end vertex. A matching is a set of mutually independent edges. 

We say a graph $G$ is \emph{good} if there exists a $\chi$-coloring of $G$ which can be reached from any coloring in $R_{\ell}(G)$ by recoloring each vertex at most $n$ times, for all $\ell \geq \chi(G)$+1. We call such a $\chi$-coloring a \emph{good coloring} of $G$. Note that every graph with at most three vertices is good. When we speak about a path between two $\ell$-colorings $\alpha$ and $\beta$ of $G$, we mean a path between the two vertices in $R_{\ell}(G)$ corresponding to the colorings $\alpha$ and $\beta$.

\section{Properties of recolorable graphs}

We first make the following observation.
\begin{proposition}\label{good}
    If a graph $G$ is good, then $G$ is recolorable with $\ell$ recoloring diameter at most 2$n^2$ for all $\ell \geq \chi(G)$+1.
\end{proposition}

A class of graphs $\mathcal{G}$ is called hereditary if for every $G \in \mathcal{G}$, every induced subgraph $H$ of $G$ is also in $\mathcal{G}$. We prove some results about a minimal graph $G$ which is not recolorable. Our results are mainly structured toward hereditary classes of graphs. If $G\in \mathcal{G}$ is a minimal graph that is not recolorable, then $G$ can not have any of the following properties:
\begin{itemize}
\setlength\itemsep{0em}
\item[(i)] $G$ has a vertex of degree at most $\chi(G)$-1
\item[(ii)] $G$ has two non-adjacent vertices $u$ and $v$, such that $N(u)\subseteq N(v)$
\item[(iii)] $G$ is the disjoint union or join of two graphs
\item[(iv)] $G$ has a tight clique cutset \cite{Belavadi}
\end{itemize}

We prove each of these properties below. In \cite{biedl2021}, some of these properties were applied, starting with a single vertex, to build a class of graphs known as OAT-graphs that are recolorable. If $G$ is an OAT-graph, then the $\ell$ recoloring diameter of $G$ is at most 4$n^2$ for all $\ell\geq \chi(G)$+1 \cite{biedl2021}. Our approach not only enables us to characterize graphs that are recolorable but also graphs that are not recolorable. From now on we assume $\ell \geq \chi(G)$+1. First we shall see a result on isomorphic $\chi$-colorings.

\begin{lemma}[Renaming Lemma \cite{bonamy2018}]
\label{lem:recolor}
Let $\alpha'$ and $\beta'$ be two $\chi$-colorings of $G$ that induce the same partition of vertices into color classes and let $\ell \ge \chi(G)+1$. Then $\alpha'$ can be recolored into $\beta'$ in $R_{\ell}(G)$ by recoloring each vertex at most twice.
\end{lemma}

 The Renaming Lemma implies that we can recolor any coloring $\alpha$ of a complete graph to a coloring $\beta$ by recoloring each vertex at most twice. Therefore, every complete graph is good.

\begin{lemma}\label{degree}
Let $G$ contain a vertex $v$ of degree at most $\chi(G)$-1. Then, either\\
(i) ~$G$ is recolorable or\\
(ii) $G$-$v$ is not recolorable.
\end{lemma}
\begin{proof}
Let $G$ be a graph and let $v$ be a vertex in $G$ with $d(v)\le \chi(G)$-1. Assume that $G$-$v$ is recolorable. Let $\ell\geq \chi(G)$+1. Let $\alpha$ and $\beta$ be any two $\ell$-colorings of $G$. We prove that there exists a path between them in $R_{\ell}(G)$. Let $\alpha_1$ and $\beta_1$ denote the restriction of $\alpha$ and $\beta$ to $G$-$v$, respectively. Since $G$-$v$ is recolorable, there exists a path between $\alpha_1$ and $\beta_1$, say $P$, in $R_{\ell}(G$-$v)$. Since $\ell \geq \chi(G)+1 \geq d(v)$+2, for each coloring in $R_{\ell}(G)$, there is at least one color that does not appear in $N[v]$ in that coloring. Thus every $\ell$-coloring of $G$-$v$ can be extended to an $\ell$-coloring of $G$.\\
\textit{Claim}: There exists a path between $\alpha$ and $\beta$ in $R_{\ell}(G)$.\\
\\
\textit{Proof of Claim}: Let $\gamma$ and $\delta$ be any two $\ell$-colorings of $G$, such that the restriction of $\gamma$ to $G$-$v$, say $\gamma_1$, and the restriction of $\delta$ to $G$-$v$, say $\delta_1$, are adjacent in the path $P$. We prove the claim by showing that there exists a path between $\gamma$ and $\delta$ in $R_{\ell}(G)$. Since $\gamma_1$ and $\delta_1$ are adjacent in $P$, they differ on exactly one vertex, say $u$, in $G$-$v$.\\
\textit{Case 1}: If $u\notin N(v)$ or $\delta_{1}(u) \neq \gamma(v)$, then starting with the coloring $\gamma$, recolor vertex $u$ with color $\delta(u)$ and recolor vertex $v$ with color $\delta(v)$ to obtain the coloring $\delta$.\\
\textit{Case 2}: If $u\in N(v)$ and $\delta_{1}(u)$ = $\gamma(v)$, then since $\ell \geq d(v)$+2, there exists a color, say $r$, that does not appear in $\gamma(N(v))\cup \delta(N(v))$. Hence starting with the coloring $\gamma$, recolor vertex $v$ with color $r$, recolor $u$ with color $\delta(u)$, and recolor $v$ with color $\delta(v)$ to obtain the coloring $\delta$. 
\end{proof}

Since every graph $G$ is ($\Delta(G)$+1)-colorable, the above result shows that every graph $G$ is ($\Delta(G)$+2)-mixing, which was first proved in \cite{Jerrum1995}.

\begin{lemma}\label{independent}
Let $G$ be a graph with non-adjacent vertices $u$ and $v$, such that $N(u)\subseteq N(v)$. Then either\\
(i) ~$G$ is recolorable or\\
(ii) $G$-$u$ is not recolorable.\\
Furthermore, if $G$-$u$ is good then $G$ is good.
\end{lemma}
\begin{proof}
Let $G$ be a graph that contains two non-adjacent vertices $u$ and $v$, such that $N(u)\subseteq N(v)$ and assume that $G$-$u$ is recolorable. Let $\ell\geq \chi(G)$+1 and let $\alpha$ be any $\ell$-coloring of $G$. Let $\beta$ be a $\chi$-coloring of $G$ such that $\beta(u)$ = $\beta(v)$. We prove that there is path between $\alpha$ and $\beta$ in $R_{\ell}(G)$ and hence $G$ is recolorable. Let $\alpha_1$ and $\beta_1$ be restrictions of $\alpha$ and $\beta$, respectively, to $G$-$u$. Since $G$-$u$ is recolorable, there exists a shortest path, say $P$, between $\alpha_1$ and $\beta_1$ in $R_{\ell}(G$-$u)$. Note that every coloring in the path $P$ can be extended to a coloring of $G$ by coloring $u$ the color of $v$. Let $P^*$ be a path in $R_{\ell}(G)$ that is described as follows; start with the coloring $\alpha$, recolor vertex $u$ with the color of $v$ and then follow the recoloring sequence in $P$ while recoloring $u$ with the color of $v$ whenever $v$ is recolored in $P$. Hence $P^{*}$ is a path between $\alpha$ and $\beta$ in $R_{\ell}(G)$.

We now prove that if $G$-$u$ is good, then $G$ is good. We retroactively choose $\beta_1$ to be a good coloring of $G$-$u$. Then in the path $P$ every vertex of $G$-$u$ is recolored at most $n$-1 times, so in the path $P^{*}$ every vertex is recolored at most $n$ times. Therefore $\beta$ is a good coloring of $G$.
\end{proof}

Before we move further, we make the following remark.\\
\\
\noindent
\textbf{Remark 1}: Let $G$ = ($V$, $E$) be a graph and let $\alpha$ be a coloring of $G$. Let $S$ be a set of colors such that $\alpha(V)\subseteq S$. Let $\alpha^*$ be a coloring of $G$ such that there is a path $P$ from $\alpha$ to $\alpha^*$ where every coloring in $P$ only uses the colors in $S$. Let $S^*$ be another set of colors such that $|S|$ = $|S^*|$. Define a bijection $f$ : $S\setminus S^* \to S^*\setminus S$. Let $\beta$ and $\beta^*$ be two colorings of $G$ such that, for all $v\in V(G)$, 
\begin{equation*}
    \beta(v) =
  \begin{cases}
                                   \alpha(v) & \textit{if } \alpha(v) \in S\cap S^* \\
                                   f(\alpha(v)) & \textit{if } \alpha(v) \notin S\cap S^* ;
  \end{cases}
  \hspace{10mm}
  \beta^*(v) =
  \begin{cases}
                                   \alpha^*(v) & \textit{if } \alpha^*(v) \in S\cap S^* \\
                                   f(\alpha^*(v)) & \textit{if } \alpha^*(v) \notin S\cap S^* ;
  \end{cases}
\end{equation*}
We show that there is a path $P^*$ from $\beta$ to $\beta^*$ where every coloring in $P^{*}$ only uses the colors in $S^*$ as follows. Whenever we recolor a vertex $v$ with color $c$ in the path $P$, $P^*$ is obtained by recoloring the vertex $v$ with color $c$ if $c\in S\cap S^*$, or with color $f(c)$ if $c\notin S\cap S^*$. Furthermore, for each vertex, the number of times it is recolored in path $P$ and the number of times it is recolored in path $P^*$ are the same.

\begin{lemma}\label{join}
Let $G$ be the disjoint union or join of two graphs $G_1$ and $G_2$, then either\\
(i) ~$G$ is recolorable or\\
(ii) $G_i$ is not recolorable for some $i \in \{1, 2\}$.\\
Furthermore, if $G_1$ and $G_2$ are good then $G$ is good.
\end{lemma}
\begin{proof}
Let $G$ be the disjoint union or join of two graphs $G_1$ and $G_2$. Assume that $G_i$ is recolorable for all $i \in \{1, 2\}$. Let $\ell\geq \chi(G)$+1 and let $S$ = $\{1,\dots,\ell\}$ be the set of available colors. Let $\alpha$ be any $\ell$-coloring of $G$ and let $\alpha_i$ denote the restriction of $\alpha$ to $G_i$, for $i$ = 1, 2. Let $V(G_i)$ = $V_i$ and let $n_i$ = $|V_i|$, for $i$ = 1, 2.\\
\\
\noindent
\textit{Case 1}: \textbf{$G$ is the disjoint union of $G_1$ and $G_2$}.\\
 Let $\beta_1$ be a coloring of $G_1$ with colors in $\{1,\dots,\chi(G_1)\}$ and let $\beta_2$ be a coloring of $G_2$ with colors in $\{1,\dots,\chi(G_2)\}$. Let $\beta$ be a $\chi$-coloring of $G$ such that $\beta(u)$ = $\beta_i(u)$ if $u\in V_i$, where $i\in \{1,2\}$. We show that there is a path between $\alpha$ and $\beta$ in $R_{\ell}(G)$. Since $G_i$ is recolorable, for $i$ = 1, 2, there exists a path between $\alpha_i$ and $\beta_i$ in $R_{\ell}(G_i)$. Let $P_i$ be one such path between $\alpha_i$ and $\beta_i$ in $R_{\ell}(G_i)$, for $i$ = 1, 2, where the number of times any vertex is recolored is minimum. Since $G$ is a disjoint union, the recoloring sequence of $P_1$, with the colors on $G_2$ unchanged, followed by the recoloring sequence of $P_2$, with the colors on $G_1$ unchanged, will correspond to a path, say $P$, between $\alpha$ and $\beta$ in $R_{\ell}(G)$. 

We now prove that if $G_1$ and $G_2$ are good, then $G$ is good. We retroactively choose $\beta_1$ to be a good coloring of $G_1$ with colors 1,$\dots,\chi(G_1)$ and choose $\beta_2$ to be a good coloring of $G_2$ with colors 1,$\dots,\chi(G_2)$. We prove that $\beta$ is a good coloring of $G$. Since $\beta_i$ is a good coloring of $G_i$, in $P_i$ each vertex of $G_i$ is recolored at most $n_i$ times, for $i$ = 1, 2. Therefore, in the path $P$ each vertex of $G$ is recolored at most $max\{n_1, n_2\}< n$ times.\\
\\
\noindent
\textit{Case 2}: \textbf{$G$ is the join of $G_1$ and $G_2$}.\\
Note that $\chi(G)$ = $\chi(G_1)$ + $\chi(G_2)$. Let $\beta_1$ be a $\chi$-coloring of $G_1$ with colors in $\{1,\dots,\chi(G_1)\}$ and let $\beta_2$ be a $\chi$-coloring of $G_2$ with colors in $\{\chi(G_1)$+1$,\dots,\chi(G)\}$. Let $\beta$ be the $\chi$-coloring of $G$ such that $\beta(u)$ = $\beta_i(u)$ if $u\in V_i$, where $i\in \{1,2\}$. We show that there is a path between $\alpha$ and $\beta$ in $R_{\ell}(G)$. In any $\ell$-coloring of $G$, the colors that appear on the vertices of $G_1$ and $G_2$ are mutually distinct. So in any $\ell$-coloring of $G$, for some $i\in \{1,2\}$, at least $\chi(G_i)$+1 colors are available to color the vertices of $G_i$. Without loss of generality, assume that there are at least $\chi(G_1)$+1 colors available to color the vertices of $G_1$ under $\alpha$. Since $G_1$ is recolorable, by Remark 1, we can recolor $\alpha_1$ to a $\chi$-coloring of $G_1$, say $\beta_{1}^{*}$, isomorphic to $\beta_1$, using the colors that do not appear on the vertices of $G_2$. Let $P_1$ be one such path, between $\alpha_1$ and $\beta_1^{*}$, where the number of times any vertex is recolored is minimum. Now we have at least $\chi(G_2)$+1 colors that do not appear on the vertices of $G_1$. Since $G_2$ is recolorable, use these colors to recolor the vertices of $G_2$ from coloring $\alpha_2$ to a $\chi$-coloring of $G_2$, say $\beta_2^*$, isomorphic to $\beta_2$ using the colors that do not appear on the vertices of $G_1$. Let $P_2$ be one such path, between $\alpha_2$ and $\beta_2^{*}$, where the number of times any vertex is recolored is minimum. Since $G$ is the join of $G_1$ and $G_2$, the recoloring sequence of $P_1$, with the colors on $G_2$ unchanged, followed by the recoloring sequence of $P_2$, with the colors on $G_1$ unchanged, will correspond to a path, say $P$, in $R_{\ell}(G)$ between $\alpha$ and a $\chi(G)$-coloring $\beta^{*}$ of $G$ which is isomorphic to $\beta$. \\
\\
Let $n_i$ = 1, for some $i\in \{1,2\}$. Without loss of generality, let $V(G_2)$ = $\{v\}$. By Remark 1, we can assume there is at most one different color used by $\beta_1^*$ compared to colors used by $\beta_1$. Let $c^* \in \beta_1^*(V_1)\setminus \beta_1(V_1)$ and $c \in \beta_1(V_1)\setminus \beta_1^*(V_1)$. Since $\ell \geq \chi(G_1)$+2, there exist some color, say $d$, not in $\beta_1(V_1)\cup \beta_1^*(V_1)$. Now recolor the vertex $v$ with color $d$, recolor every vertex in $G_1$ colored $c^*$ with color $c$, and recolor $v$ with color $\chi(G_1)$+1 to obtain the coloring $\beta$. Therefore there is a path between $\alpha$ and $\beta$ in $R_{\ell}(G)$. We now prove that if $G_1$ is good, then $G$ is good. We retroactively choose $\beta_1$ to be a good coloring of $G_1$ that uses the colors in $\{1,\dots,\chi(G_1)\}$. Let $\beta$ be the $\chi$-coloring of $G$ obtained by extending the coloring $\beta_1$ of $G_1$ by coloring the vertex in $G_2$ the color $\chi(G_1)$+1. We prove that $\beta$ is a good coloring of $G$. Since $\beta_1$ is a good coloring of $G_1$, it follows from Remark 1 that in the path $P_1$ every vertex was recolored at most $n_1$ times. Therefore we can recolor $\alpha$ to $\beta$ by recoloring every vertex at most $max\{n_1$+1, 2$\} \leq n$ times.\\
\\
Let $n_i > 1$, for all $i\in  \{1,2\}$. By the Renaming Lemma, we can reach $\beta$ from $\beta^{*}$ by recoloring each vertex at most twice. Thus there is a path between $\alpha$ and $\beta$ in $R_{\ell}(G)$. We now prove that if $G_1$ and $G_2$ are good, then $G$ is good. We retroactively choose $\beta_1$ to be a good coloring of $G_1$ with colors in $\{1,\dots,\chi(G_1)\}$ and choose $\beta_2$ to be isomorphic to a good coloring of $G_2$ with colors in $\{\chi(G_1)$+1,$\dots,\chi(G)\}$. Since $\beta_1$ is a good coloring of $G_1$, it follows from Remark 1 that we can recolor $\alpha_1$ to $\beta_{1}^{*}$ by recoloring each vertex of $G_1$ at most $n_1$ times. Since $\beta_2$ is isomorphic to a good coloring, it follows from Remark 1 that we can recolor $\alpha_2$ to $\beta_2^*$ by recoloring every vertex of $G_2$ at most $n_2$ times. We can apply the Renaming Lemma by recoloring each vertex of $G$ at most twice. Therefore we can reach $\beta$ from $\alpha$ in $R_{\ell}(G)$ by recoloring each vertex at most $max\{n_1, n_2\}$+2 $\leq n$ times.\\
 \end{proof}

\section{Triangle-free graphs}\label{sec:traingle}

Every $P_4$-free graph $G$ with at least two vertices is either a join or a disjoint union of some proper induced subgraphs $G_1$ and $G_2$ of $G$. The proof of Theorem 1 in \cite{biedl2021} can be altered to obtain the following result.

\begin{lemma}(\cite{biedl2021})\label{lem1}
    Let $G$ be a $P_4$-free graph. Then we can recolor any $\ell$-coloring $\alpha$ of $G$ to any $\ell$-coloring $\beta$ of $G$ by recoloring each vertex at most 4 times, for all $\ell\geq \chi(G)$+1.
\end{lemma}

The above result proves that every $P_4$-free graph is recolorable with $\ell$ recoloring diameter at most 4$n$. Since complete multipartite graphs are precisely the ($P_2$+$P_1$)-free graphs, which is a subclass of $P_4$-free graphs, we have the following.


\begin{theorem}\label{multipartite}
Every complete multipartite graph $G$ is good and hence recolorable. Furthermore, the $\ell$ recoloring diameter of $G$ is at most 4n for all $\ell\geq \chi(G)$+1.
\end{theorem}
\begin{proof}
It is easy to see that the result holds for graphs with at most 3 vertices. For graphs with at least 4 vertices, the result follows from Lemma \ref{lem1}.
 \end{proof}

\begin{theorem}(\cite{olariu})\label{paw-triangle}
A connected graph $G$ is paw-free if and only if it is triangle-free or a complete multipartite graph.
\end{theorem}

From Theorems \ref{multipartite} and \ref{paw-triangle}, we have the following theorem.
\begin{theorem}
Every (triangle, $H$)-free graph is recolorable if and only if every (paw, $H$)-free graph is recolorable. Furthermore, every (triangle, $H$)-free graph is good if and only if every (paw, $H$)-free graph is good.
\end{theorem}

In \cite{bonamy2014}, it was proved that the diameter of the reconfiguration graph of 3-colorings of $P_n$, denoted by $R_3(P_n)$, is $\Omega(n^2)$. We prove a result on cycles.

\begin{lemma}\label{cycle}
    For all $p\geq 4$, the graph $C_n$, is $p$-mixing with $p$ recoloring diameter at most 4$n$.
\end{lemma}
\begin{proof}
    Let $C_n$ denote a cycle on vertices $v_1,\dots,v_n$ with edges $v_iv_{i+1}$ where subscripts are mod $n$. Let $\gamma$ be a 3-coloring of $G$ coloring vertices $v_1,\dots,v_{n-1}$ with colors 1, 2, 3, 1, \dots and by coloring $v_n$ with color $n$ (mod 3) if $n\not\equiv$ 1 (mod 3) or with color 2 otherwise. We prove that any $p$-coloring $\alpha$ of $C_n$, $p\geq 4$, can be recolored to $\gamma$ by recoloring each vertex at most twice.

    Starting with coloring $\alpha$, recolor vertex $v_1$ with color 1 if $\alpha(v_n) \neq$ 1 $\neq \alpha(v_2)$. If $v_j$ is colored 1 for $j\in \{n, 2\}$, then since $v_{j}$ has only two neighbors and at least four available colors, there exists a color $c$ that does not appear in its closed neighborhood. Recolor $v_{j}$ with color $c$ and recolor vertex $v_1$ with color 1. For every vertex in ($v_2,\ v_3,\dots, v_{n-1}$) in the order of vertices: Recolor the vertex $v_i$ with color $i$ (mod 3) if $\alpha(v_{i+1}) \neq i$ (mod 3). If $\alpha(v_{i+1})$ = $i$ (mod 3), then since $v_{i+1}$ has only two neighbors and at least four available colors, there exists a color $c$ that does not appear in its closed neighborhood. Recolor $v_{i+1}$ with color $c$ and recolor vertex $v_i$ with color $i$ (mod 3). Now recolor vertex $v_n$ with color $n$ (mod 3) if $n\not\equiv$ 1 (mod 3) or else recolor it with color 2. Thus we can recolor $\alpha$ to the coloring $\gamma$ by recoloring each vertex at most twice.

    Thus, given any two $p$-colorings $\alpha$ and $\beta$ of $C_n$, there exists a path from $\alpha$ to $\gamma$ and a path from $\beta$ to $\gamma$, each of length at most 2$n$. Therefore there exists a path between $\alpha$ and $\beta$, in $R_{p}(G)$, of length at most 4$n$.
\end{proof}

\begin{theorem}
Every component of a $(triangle, claw)$-free graph $G$ is recolorable with $\ell$ recoloring diameter at most 2$n^{2}$ for all $\ell\geq \chi(G)$+1 or isomorphic to $C_{2q}$ for some $q\geq 3$. 
\end{theorem}
\begin{proof}
Let $G$ be a connected (triangle, claw)-free graph. Then $G$ is an induced path or a cycle \cite{hamel2019}. If $G$ is isomorphic to $C_{2q}$ for some $q\geq 3$, then $R_{3}(G)$ is disconnected \cite{Cereceda2008} and hence $G$ is not recolorable. If $G$ is an odd cycle, then $\chi(G)$ = 3 and the result follows from Lemma \ref{cycle}. If $G$ is an induced path on at least three vertices or a 4-cycle, then we prove that $G$ is good by induction on $n$. Path $G$ contains two independent vertices $u$ and $v$ such that $N(u)\subseteq N(v)$ and $G$-$u$ is a path. By the induction hypothesis, $G$-$u$ is good and the result follows from Lemma \ref{independent}. 
\end{proof}

\begin{corollary}
Every component of a $(paw, claw)$-free graph $G$ is recolorable with $\ell$ recoloring diameter at most 2$n^{2}$ for all $\ell\geq \chi(G)$+1 or isomorphic to $C_{2q}$ for some $q\geq 3$. 
\end{corollary}

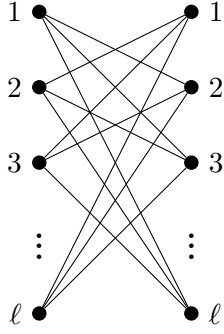
\begin{figure}[h]
\centering
\begin{tikzpicture}[scale=0.5]
\tikzstyle{vertex}=[circle, draw, fill=black, inner sep=0pt, minimum size=5pt]

    \node[vertex, label=left:1](1) at (0,4) {};
    \node[vertex, label=left:2](2) at (0,2) {};
    \node[vertex, label=left:3](3) at (0,0) {};
	\node[vertex, label=left:$\ell$](4) at (0,-4) {};
	\node[vertex, label=right:1](5) at (4,4) {};
	\node[vertex, label=right:2](6) at (4,2) {};
	\node[vertex, label=right:3](7) at (4,0) {};
	\node[vertex, label=right:$\ell$](8) at (4,-4) {};
	
	\node at (0,-2) {\textbf{\vdots}};
	\node at (4,-2) {\textbf{\vdots}};
     
    \draw(1)--(6);
    \draw(1)--(7);
    \draw(1)--(8);
    \draw(2)--(5);
    \draw(2)--(7);
    \draw(2)--(8);
    \draw(3)--(5);
    \draw(3)--(6);
    \draw(3)--(8);
    \draw(4)--(5);
    \draw(4)--(6);
    \draw(4)--(7);
    
\end{tikzpicture}
\caption{A frozen $\ell$-coloring of the graph $K_{\ell, \ell}$-$M$. \cite{Cereceda2008}}
\label{fig:frozenbipartite}
\end{figure}

\begin{theorem}
Every component of a $(triangle, co$-$diamond)$-free graph $G$ is either recolorable with $\ell$ recoloring diameter at most 6$n$ for all $\ell\geq \chi(G)$+1 or isomorphic to $K_{\ell, \ell}$-$M$, where $M$ is a maximum matching for some $\ell\geq \chi(G)$+1.
\end{theorem}
\begin{proof}
Let $G$ be a connected (triangle, co-diamond)-free graph. There are two cases:\\
\\
\textit{Case 1}: \textbf{$G$ is a bipartite graph}.\\
\noindent
If $G$ is isomorphic to $K_{\ell, \ell}$-$M$, for some $\ell \geq \chi(G)$+1, then $G$ has a frozen $\ell$-coloring \cite{Cereceda2008}. If $G$ is not isomorphic to $K_{\ell, \ell}$-$M$, for any $\ell \geq \chi(G)$+1, then we prove that there is a path between any two $\ell$-colorings $\alpha$ and $\beta$ of $G$ in $R_{\ell}(G)$ where each vertex is recolored at most 6 times. Let $A$ and $B$ be the partite sets of $V(G)$. Let $C$ = $\{1,\dots,\ell\}$ be the set of available colors.

Consider an $\ell$-coloring $\alpha$ of $G$. Let us assume there is a color $i\in [\ell]$ that does not appear on any of the vertices in one of the partite sets, say $A$, under $\alpha$. Then we recolor each vertex in $B$ with the color $i$. Next recolor each vertex in $A$ with one of the available colors $j \neq i$ to obtain a 2-coloring $\alpha_1$ of $G$. Thus there is a path between $\alpha$ and $\alpha_1$ in $R_{\ell}(G)$ where each vertex is recolored at most once. Now assume that all available colors appear in both sets $A$ and $B$ under $\alpha$. For each $i\in [\ell]$, let $A_i$ = $\{u \in A| \alpha(u)$ = $i\}$; then $A_1$, $A_2,\dots,\ A_{\ell}$ is a partition of $A$ into non-empty sets. Define $B_1$, $B_2$,\dots, $B_{\ell}$ similarly. Since the vertices in $A_i$ and $B_i$ receive the same color under $\alpha$, $A_{i}$ is anticomplete to $B_i$ for all $i\in [\ell]$. If $|A_i| > 1$, for some $i\in [\ell]$, then let $u$, $v \in A_{i}$ and $b_i\in B_i$. Then for some $x$, a neighbor of $b_i$ in $A$, $\{b_i, x, u, v\}$ induces a co-diamond, a contradiction. Thus, for all $i\in [\ell]$, $|A_i|$ = 1, and similarly $|B_i|$ = 1. If for all $i$, $j\in [\ell]$, $i\neq j$, $A_i$ is complete to $B_j$, then $G$ is isomorphic to $K_{\ell, \ell}$-$M$. If $A_i$ is anticomplete to some $B_j$, $i\neq j$, then by recoloring the vertex in $A_i$ with the color $j$ we obtain a coloring in which the color $i$ does not appear on any vertex in $A$. Then we can use the procedure explained above to obtain a 2-coloring $\alpha_1$ of $G$. Thus there is a path between $\alpha$ and $\alpha_1$ in $R_{\ell}(G)$ where each vertex is recolored at most twice.

Similarly, if $\beta$ is an $\ell$-coloring of $G$ then there is a path in $R_{\ell}(G)$ between $\beta$ and a 2-coloring $\beta_1$ of $G$ where each vertex is recolored at most twice. Since $\alpha_1$ and $\beta_1$ are isomorphic 2-colorings of $G$, by the Renaming Lemma, there is a path between them in $R_{\ell}(G)$ where each vertex is recolored at most twice. This proves that there is a path between $\alpha$ and $\beta$ in $R_{\ell}(G)$ of length at most 6$n$.\\
\\
\textit{Case 2}: \textbf{$G$ is not a bipartite graph}. \\
\noindent
Then $G$ has an induced odd cycle with at least five vertices. Since $G$ is co-diamond-free, $G$ is $C_p$-free for all $p\geq 7$. Thus $G$ contains an induced 5-cycle, say $C$, with $V(C)$ = $\{v_1,\dots,v_5\}$ such that $v_iv_{i+1}\in E(G)$, for all $i$ (mod 5). If $G$ is not isomorphic to $C_5$, then since $G$ is connected there is a vertex $u_1$ not in $C$ adjacent to some vertex in $C$. Without loss of generality, we assume that $u_1$ is adjacent to $v_1$. Since $G$ is triangle-free, $\{v_5, v_2, u_1\}$ must be an independent set. If $u_1$ is not adjacent to $v_4$, then $\{u_1, v_2, v_4, v_5\}$ induces a co-diamond, a contradiction. So $u_1$ must be adjacent to $v_4$, but then $\{u_1, v_5, v_2, v_3\}$ induces a co-diamond, a contradiction. Therefore, $G$ is isomorphic to $C_5$ and the result follows from Lemma \ref{cycle}. 
\end{proof}

\begin{corollary}
Every component of a $(paw, co$-$diamond)$-free graph $G$ is either recolorable with $\ell$ recoloring diameter at most 6$n$ for all $\ell\geq \chi(G)$+1 or isomorphic to $K_{\ell, \ell}$-$M$, where $M$ is a maximum matching for some $\ell\geq \chi(G)$+1.
\end{corollary}

Note that we can check if a graph $H$ is isomorphic to some $K_{\ell, \ell}$-$M$ in polynomial time. Thus we can check in polynomial time, for a (triangle, co-diamond)-free graph $G$, whether $R_{\ell}(G)$ is connected.

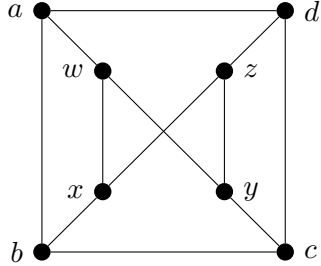
\begin{figure}[h]
\centering
\begin{tikzpicture}[scale = 0.4,
every edge/.style = {draw=black,very thick},
 vrtx/.style args = {#1/#2}{%
      circle, draw, thick, fill=black, inner sep = 0pt,
      minimum size=2mm, label=#1:#2}
                    ]

\node(1) [vrtx=right/$z$] at (2,2) {};
\node(2)[vrtx=right/$y$] at (2,-2) {};
\node(3)[vrtx=left/$x$] at (-2,-2) {};
\node(4)[vrtx=left/$w$] at (-2,2) {};
\node(5)[vrtx=right/$d$] at (4,4) {};
\node(6)[vrtx=right/$c$] at (4,-4) {};
\node(7)[vrtx=left/$b$] at (-4,-4) {};
\node(8)[vrtx=left/$a$] at (-4,4) {};
   \draw(1)--(2);
    \draw(1)--(3);
    \draw(1)--(5);
    \draw(2)--(4);
    \draw(2)--(6);
    \draw(3)--(4);
    \draw(3)--(7);
    \draw(4)--(8);
    \draw(5)--(6);
    \draw(5)--(8);
    \draw(6)--(7);
    \draw(7)--(8);
\end{tikzpicture}
\caption{The graph $F$.} 
\label{fig:graphF}
\end{figure}

\begin{lemma}\label{strF}
Let $G$ be a connected 3-regular $(triangle, 4K_1)$-free graph that contains an induced $P_3$+$P_1$. Then $G$ has two independent vertices $u$, $v$ such that $N(u)\subseteq N(v)$ or $G$ is isomorphic to the graph $F$ of Figure \ref{fig:graphF}. 
\end{lemma}
\begin{proof}
Let $G$ be a graph that satisfies the hypothesis and let $\{v_1, a_{12}, v_2, v_3\}$ induce a $P_3$+$P_1$ in $G$ such that $\{v_1a_{12}, a_{12}v_2\}\subseteq E(G)$. For all $i$ (mod 3), define\\
$A_{i\ i+1}$ = $\{v \in V(G)\mid v$ is complete to $\{v_i, v_{i+1}\}$ and not adjacent to $v_{i+2}\}$\\
$B_{i}$ = $\{v\in V(G)\mid v$ is adjacent to $v_i$ and anticomplete to $\{v_{i+1}, v_{i+2}\}\}$, and\\
$D$ = $\{v\in V(G)\mid v$ is complete to $\{v_1, v_2, v_3\}\}$.\\
All subscripts are mod 3. Since $G$ is $4K_1$-free, every vertex in $V(G)\setminus \{v_1, v_2, v_3\}$ has a neighbor in $\{v_1, v_2, v_3\}$ and hence is in $A_{i\ i+1}$ or $B_i$ or $D$ for some $i \in [3]$. Since $G$ is triangle-free, for all $i\in [3]$, we see that
\begin{itemize}
\setlength\itemsep{0em}
\item[(i)] $A_{12}\cup A_{23} \cup A_{31} \cup D$ is an independent set,
\item[(ii)] $B_i$ is an independent set,
\item[(iii)] $A_{i\ i+1}$ is anticomplete to $B_i \cup B_{i+1}$, and
\item[(iv)] $B_i$ is anticomplete to $D$.
\end{itemize}

If for some $i\in [3]$, $b_i\in B_i$ has two neighbors $x$ and $y$ in $A_{i+1\ i+2}$, then since $G$ is 3-regular, we have two independent vertices $x$ and $y$ such that $N(x)$ = $N(y)$ = $\{v_{i+1}, v_{i+2}, b_i\}$ and the result follows. So we assume that every vertex in $B_i$ has at most one neighbor in $A_{i+1\ i+2}$. If for some $i\in [3]$, $B_{i}$ contains more than one vertex, say $u$ and $v$, then $\{u, v, v_{i+1}, v_{i+2}\}$ induces a $4K_1$, a contradiction. Therefore, for all $i\in [3]$, $B_i$ contains at most one vertex. \\

By construction, we have $a_{12}\in A_{12}$. Since $d(a_{12})$ = 3, $a_{12}$ must have a neighbor in $B_3$, say $b_3$. Since $d(b_3)$ = 3 and $b_3$ can have at most one neighbor in $A_{12}$, $b_3$ must have a neighbor in $B_1 \cup B_2$. Without loss of generality, let $b_1\in B_1$ be adjacent to $b_3$. We have two cases.\\
\\
\textit{Case 1}: \textbf{$A_{31} \neq \emptyset$}. \\
\noindent
Let $a_{31} \in A_{31}$. Then $N(v_1)$ = $\{b_1, a_{12}, a_{31}\}$. Hence $A_{12}$ = $\{a_{12}\}$, $A_{31}$ = $\{a_{31}\}$, and $D$ = $\emptyset$. Since $d(a_{31})$ = 3, $a_{31}$ must be adjacent to some vertex $b_2 \in B_2$. If $b_1$ and $b_2$ are not adjacent, then $\{b_1, b_2, v_3, a_{12}\}$ induces a $4K_1$, a contradiction. Hence $N(b_1)$ = $\{v_1, b_3, b_2\}$ and $N(b_2)$ = $\{v_2, b_1, a_{31}\}$. If $a_{23}\in A_{23}$, then since $|B_1|$ = 1 and $b_1$ is not adjacent to $a_{23}$, $d(a_{23})$ = 2, a contradiction. If $A_{23}$ = $\emptyset$, then $N(v_2)$ = $\{b_2, a_{12}\}$. This implies $d(v_2) = 2$, a contradiction.\\ 
\\
\textit{Case 2}: \textbf{$A_{31}$ = $\emptyset$}. \\
\noindent
Since $d(v_1)$ = 3 and $|B_1|$ = 1, either there is a vertex $y\in A_{12}$ where $y\neq a_{12}$ or a vertex $d\in D$. If $y\in A_{12}$, then $\{b_1, a_{12}, y, v_3\}$ induces a $4K_1$, a contradiction. Therefore there is a vertex $d\in D$ and $N(v_1)$ = $\{b_1, a_{12}, d\}$. This also implies $D$ = $\{d\}$. Since $d(v_3)$ = 3 and $|B_3|$ = 1, $v_3$ must be adjacent to some vertex $a_{23} \in A_{23}$. Thus $N(v_3)$ = $\{b_3, d, a_{23}\}$. Since $d(a_{23})$ = 3 and $|B_1|$ = 1, $a_{23}$ must be adjacent to $b_1$ and hence $N(a_{23})$ = $\{v_2, v_3, b_1\}$. Therefore the graph $G$ is isomorphic to the graph $F$.
\end{proof}

\begin{lemma}\label{recolorF}
The graph $F$ is recolorable.
\end{lemma}
\begin{proof}
Let $F$ be the graph of Figure \ref{fig:graphF}. Note that $\chi(F)$ = 3. Given any $\ell$-coloring $\alpha$ of $F$, where $\ell \geq 4$, we prove that there is a 3-coloring $\beta$ of $F$ that induces the color classes $\{a,x,c\}$, $\{b,y,d\}$, and $\{w,z\}$, such that there is a path between $\alpha$ and $\beta$ in $R_{\ell}(F)$. Once we prove that, by the Renaming Lemma, there exists a path between any two isomorphic $\chi$-colorings of $F$ in $R_{\ell}(F)$, so $R_{\ell}(F)$ is connected. 

Consider an $\ell$-coloring $\alpha$ of $F$. Let $\{1, 2, 3, 4\}$ be a subset of available colors. We have two cases.\\
\textit{Case 1}: \textbf{Vertices $w$ and $z$ have the same color under $\alpha$}, say 3.\\ \noindent
Without loss of generality, let $\alpha(a)$ = 1 and let $\alpha(d)$ = 2. We can recolor the vertices of $F$ starting from $\alpha$ to a 3-coloring isomorphic to $\beta$ as follows:
\begin{itemize}
\setlength\itemsep{0em}
    \item recolor $x$ with color $1$,
    \item recolor $y$ with color $2$,
    \item recolor $c$ with color $1$,
    \item recolor $b$ with color $2$.
\end{itemize}
\textit{Case 2}: \textbf{Vertices $w$ and $z$ do not have the same color under $\alpha$}. \\ \noindent
If $\alpha(d)\neq \alpha(w)$ or $\alpha(a) \neq \alpha(z)$, then we can recolor $z$ with the color of $w$ or recolor $w$ with the color of $z$, respectively. Then follow the steps under \textit{Case 1} to obtain a coloring $\beta$. So we assume that $\alpha(d)$ = $\alpha(w)$ and $\alpha(a)$ = $\alpha(z)$. Without loss of generality, let $\alpha(z)$ = $\alpha(a)$ = 1 and $\alpha(w)$ = $\alpha(d)$ = 2. We can recolor the vertices of $F$ starting from $\alpha$ to a 3-coloring isomorphic to $\beta$ as follows:
\begin{itemize}
\setlength\itemsep{0em}
    \item recolor $b$ with color 2,
    \item recolor $c$ with color 1,
    \item recolor $x$ with color 4,
    \item recolor $y$ with color 4,
    \item recolor $z$ with color 3,
    \item recolor $w$ with color 3,
    \item recolor $y$ with color 2,
    \item recolor $x$ with color 1.
\end{itemize}
\end{proof}

In \cite{Belavadi}, it was proved that every ($P_3$+$P_1$)-free graph $G$ is recolorable with $\ell$ recoloring diameter at most 6$n$. We can infer the result below from the proof.

\begin{lemma}(\cite{Belavadi})\label{p3+p1}
    Let $G$ be a ($P_3$+$P_1$)-free graph. There is a $\chi$-coloring, say $\alpha$, of $G$ such that for any $\ell$-coloring $\beta$ of $G$, there is a path between $\alpha$ and $\beta$ in $R_{\ell}(G)$ where each vertex is recolored at most 4 times, for all $\ell\geq \chi(G)$+1.
\end{lemma}

It follows from Lemma \ref{p3+p1} that every ($P_3$+$P_1$)-free graph is good.

\begin{theorem}
Every $(triangle, 4K_1)$-free graph $G$ is recolorable or isomorphic to $C_6$.
\end{theorem}
\begin{proof}
Let $G$ be a (triangle, 4$K_1$)-free graph. If $G$ is disconnected, then every component of $G$ is 3$K_1$-free and hence every component is recolorable \cite{3k1free}. Then the result follows from Lemma \ref{join}. So we assume that $G$ is connected. Since every ($P_3$+$P_1$)-free graph is recolorable, from Lemma \ref{p3+p1}, we assume $G$ contains an induced $P_3$+$P_1$. 

If $G$ is bipartite, let $A$ and $B$ be the partite sets of $G$. Since $G$ is $4K_1$-free, $|A| \leq 3$ and $|B| \leq 3$. If $\delta(G)\leq 1$ or if $G$ contains independent vertices $u$ and $v$ such that $N(u)\subseteq N(v)$, then the result follows from Lemma \ref{degree} or Lemma \ref{independent}, respectively. If not, $G$ is 2-regular and we see that $G$ is isomorphic to $C_6$.

If $G$ is not bipartite, then $\chi(G)\geq 3$. If $\delta(G) \le$ 2, then the result follows from Lemma \ref{degree}. If not, then since $G$ is (triangle, $4K_1$)-free it must be 3-regular. From Lemma \ref{strF}, $G$ has two independent vertices $u$, $v$ such that $N(u)\subseteq N(v)$ or $G$ is isomorphic to the graph $F$ of Figure \ref{fig:graphF}. The result follows from Lemma \ref{independent} and Lemma \ref{recolorF}.
\end{proof}

\begin{corollary}
Every $(paw, 4K_1)$-free graph $G$ is recolorable or isomorphic to $C_6$.
\end{corollary}

\section {2$K_2$-free graphs}\label{sec:2K2}

\begin{theorem}\label{paw2k2}
Every $(2K_2, triangle)$-free graph $G$ is recolorable with $\ell$ recoloring diameter at most 2$n^{2}$, for all $\ell\geq \chi(G)$+1.
\end{theorem}
\begin{proof}
We prove that every (2$K_2$, triangle)-free graph $G$ is good. By Lemma \ref{join}, we may assume $G$ is connected. The proof is by induction on $n$. If $G$ is a bipartite graph, then since $G$ is $2K_2$-free there exist two independent vertices $u$ and $v$ such that $N(u)\subseteq N(v)$. Then by the induction hypothesis, $G$-$u$ is good, and the result follows from Lemma \ref{independent}. If $G$ is not a bipartite graph, then since $G$ is ($2K_2$, triangle)-free, $G$ contains an induced $C_5$, say $C$.

Let $V(C)$ = $\{v_1,\dots,v_5\}$ such that $v_iv_{i+1}\in E(G)$, for all $i$ (mod 5). All subscripts are mod 5. Choose maximal independent sets $A_1$, $A_2$,\dots,$A_5$ such that $v_{i} \in A_i$, and each $A_i$ is complete to $A_{i-1}\cup A_{i+1}$ and anticomplete to $A_{i+2}\cup A_{i+3}$. \\

\noindent
\textit{Claim}: For every $v \in A_i$ and $u \in N(v)$ either $u\in A_{i-1}$ or $u\in A_{i+1}$.\\
Since $u$ has a neighbour in $A_i$ and $G$ is triangle-free, $u$ must be anticomplete to $A_{i-1}\cup A_{i+1}$. If $u$ has no neighbor in $A_{i+2}\cup A_{i+3}$, then $\{v, u, v_{i+2}, v_{i+3}\}$ induces a $2K_2$, a contradiction. Up to symmetry, let $w\in A_{i+2}$ be a neighbor of $u$. Since $G$ is triangle-free, $u$ must be anticomplete to $A_{i+3}$. If there is a vertex $x$ in $A_{i+2}$ not adjacent to $u$, then $\{v_{i+3}, x, u, v\}$ induces a $2K_2$, a contradiction. Therefore $u$ is complete to $A_{i+2}$. If there is a vertex $z\in A_{i}$ not adjacent to $u$, then $\{v_{i-1}, z, u, v_{i+2}\}$ induces a $2K_2$, a contradiction. Therefore $u$ is complete to $A_i$. We have proved that $u$ is complete to $A_{i}\cup A_{i+2}$ and anticomplete to $A_{i+1}\cup A_{i+3} \cup A_{i-1}$. Thus, by maximality of $A_{i+1}$, $u$ must be in $A_{i+1}$. This completes the proof of the claim.\\

Since $G$ is connected, by the claim we conclude that $V(G)$ = $\bigcup_{i\in [5]} A_i$. If $|A_i|$ = 1, for all $i\in [5]$, then $G$ is a $C_5$, and the result follows from Lemma \ref{cycle}. If $|A_i| > 1$ for some $i\in [5]$, then there exist two independent vertices $u_1$ and $u_2$ in $A_i$ such that $N(u_1)$ = $N(u_2)$. By the induction hypothesis, $G$-$u_1$ is good, and the result follows from Lemma \ref{independent}.
\end{proof}

\begin{corollary}
Every $(2K_2, paw)$-free graph $G$ is recolorable with $\ell$ recoloring diameter at most 2$n^2$, for all $\ell\geq \chi(G)$+1.
\end{corollary}

\begin{theorem}
Every $(2K_2, claw)$-free graph $G$ is recolorable with $\ell$ recoloring diameter at most 2$n^2$, for all $\ell\geq \chi(G)$+1.
\end{theorem}
\begin{proof}
Let $G$ be a (2$K_2$, claw)-free graph. We assume $G$ is connected and has at least 4 vertices. We prove that $G$ is good. Since every ($P_3$+$P_1$)-free graph is good, by Lemma \ref{p3+p1} we may assume $G$ contains an induced $P_3$+$P_1$. The proof is by induction on $n$. Let $\{a,b,c,d\}\subset V(G)$ induce a $P_3$+$P_1$ such that $a$ and $c$ are neighbors of $b$. Since $G$ is connected, there is a vertex $x$ adjacent to $d$. For any vertex $v \in N(d)\setminus N(b)$, $v$ is complete to $\{a, c\}$; otherwise $\{a,b,v,d\}$ or $\{b,c,v,d\}$ induces a $2K_2$, a contradiction. However, this implies $\{a, c, d, v\}$ induces a claw, a contradiction. Therefore, $N(d)\subseteq N(b)$. By the induction hypothesis, $G$-$d$ is good and the result follows from Lemma \ref{independent}. 
\end{proof}

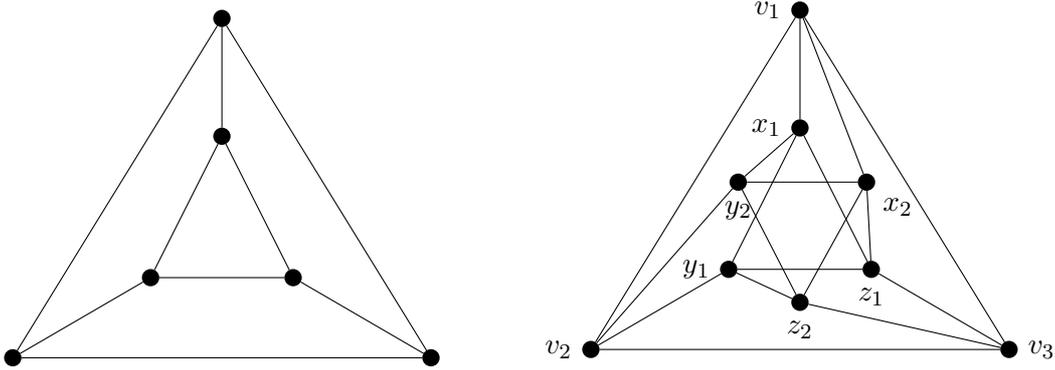
\begin{figure}[h]
\begin{center}

\begin{tikzpicture}[scale = 1.25,
every edge/.style = {draw=black,very thick},
 vrtx/.style args = {#1/#2}{%
      circle, draw, thick, fill=black, inner sep = 0pt,
      minimum size=2mm, label=#1:#2}
                    ]

\node(1)[vrtx=left/] at (2,3.75) {};
\node(2)[vrtx=left/] at (-0.2,0.15) {};
\node(3)[vrtx=left/] at (4.2,0.15) {};
\node(4)[vrtx=left/] at (2,2.5) {};
\node(5)[vrtx=left/] at (1.25,1) {};
\node(6)[vrtx=left/] at (2.75,1) {};
\draw (1) -- (2);
\draw (1) -- (3);
\draw (2) -- (3);
\draw (1) -- (4);
\draw (5) -- (2);
\draw (3) -- (6);
\draw (4) -- (5);
\draw (4) -- (6);
\draw (5) -- (6);
\end{tikzpicture} \hspace{1cm}
\begin{tikzpicture}[scale = 1.25,
every edge/.style = {draw=black,very thick},
 vrtx/.style args = {#1/#2}{%
      circle, draw, thick, fill=black, inner sep = 0pt,
      minimum size=2mm, label=#1:#2}
                    ]
\node(1)[vrtx=left/$v_1$]  at (2,3.75) {};
\node(2)[vrtx=left/$v_2$]  at (-0.2,0.15) {};
\node(3)[vrtx=right/$v_3$]  at (4.2,0.15) {};
\node(4)[vrtx=left/$x_1$]  at (2,2.5) {};
\node(5)[vrtx=left/$y_1$]  at (1.25,1) {};
\node(6)[vrtx=below/$z_1$]  at (2.75,1) {};
\node(7)[vrtx=below/$y_2$]  at (1.35,1.925) {};
\node(8)[vrtx=below/$z_2$]  at (2,0.65) {};
\node(9)[vrtx=below right/$x_2$]  at (2.7,1.925) {};
\draw (1) -- (2);
\draw (1) -- (3);
\draw (2) -- (3);
\draw (1) -- (4);
\draw (5) -- (2);
\draw (3) -- (6);
\draw (4) -- (5);
\draw (4) -- (6);
\draw (5) -- (6);
\draw (8) -- (7);
\draw (8) -- (9);
\draw (7) -- (9);
\draw (7) -- (2);
\draw (7) -- (4);
\draw (8) -- (3);
\draw (5) -- (8);
\draw (6) -- (9);
\draw (1) -- (9);

\end{tikzpicture}

\end{center}
\vspace{-2mm}
\caption{The 3-prism (left) and the 3-prism star (right).}\label{fig:3-prism}
\end{figure}

\begin{lemma}\label{prism}
    The 3-prism star graph is good and hence recolorable.
\end{lemma}
\begin{proof}
    Let $G$ be the 3-prism star graph of Figure \ref{fig:3-prism}. Note $\chi(G)$ = 3 and let $\ell \geq 4$. Let $\alpha$ be an $\ell$-coloring of $G$. Let $\beta$ be a $\chi$-coloring which induces color classes $\{v_1,z_1,z_2\}$, $\{v_2,x_1,x_2\}$, and $\{v_3,y_1,y_2\}$. We prove that there is a path from $\alpha$ to a $\chi$-coloring that is isomorphic to $\beta$ in $R_{\ell}(G)$, where each vertex is recolored at most twice. Hence by the Renaming Lemma, there is a path from $\alpha$ to $\beta$ in $R_{\ell}(G)$, where each vertex is recolored at most 4 times.

    Let $\{1,2,\dots,\ell\}$ be the set of available colors. Consider the coloring $\alpha$. Without loss of generality, let $\alpha(v_1)$ = 1, let $\alpha(v_2)$ = 2, and let $\alpha(v_3)$ = 3. Starting with the coloring $\alpha$, recolor the vertices as follows.

    If $\alpha(x_1)$ = 3, then do the following. Recolor $x_2$ with color 3, recolor $z_1$ and $z_2$ with color 2, recolor $y_1$ and $y_2$ with color 4, recolor $z_1$ and $z_2$ with color 1, recolor $x_1$ and $x_2$ with color 2, and recolor $y_1$ and $y_2$ with color 3. This is a coloring isomorphic to $\beta$ and it was obtained by recoloring each vertex at most twice.

    If $\alpha(x_1) \neq$ 3, then let $\alpha(x_1) = r \neq$ 1, for some $r\in [\ell]$. Recolor $y_1$ with color 3, recolor $z_1$ with color 1, and recolor $x_1$ with color 2. If $\alpha(x_2)$ = 3, then recolor $z_2$ with color 2, recolor $y_2$ with color 4, recolor $z_2$ with color 1, recolor $x_2$ with color 2, and recolor $y_2$ with color 3. This is a coloring isomorphic to $\beta$ and it was obtained by recoloring each vertex at most twice. If $\alpha(x_2) \neq$ 3, then recolor $y_2$ with color 3, recolor $z_2$ with color 1, and recolor $x_2$ with color 2. This is a coloring isomorphic to $\beta$ and it was obtained by recoloring each vertex at most twice.
\end{proof}

\begin{lemma}\label{clique3}
    Every ($2K_2$, diamond)-free graph $G$ with clique number 3 is recolorable with $\ell$ recoloring diameter at most 2$n^2$, for all $\ell\geq \chi(G)$+1.
\end{lemma}
\begin{proof}
We prove that every (2$K_2$, diamond)-free graph with clique number 3 is good. If possible, let $G$ be a vertex-minimal counter example. Then $G$ is not complete and, by Lemma \ref{join}, $G$ is connected. Let $G$ contain a triangle induced by the subset of vertices $Q$ = $\{v_1,v_2,v_3\}$. For all $i\in$ [3], define $B_i$ = $\{v\in V(G)\mid N(v) \cap Q$ = $\{v_i\}\}$. All subscripts are mod 3. Since $G$ is diamond-free with clique number 3, every neighbor of $v_i$ not in $Q$ is in $B_i$. \\
\textit{Claim}: $B_i\neq \emptyset$, for all $i\in$ [3].\\
\textit{Proof of Claim}: If possible, let $B_i$ = $\emptyset$ for some $i\in$ [3]. Since $G$ is connected and not complete, there is some $B_j \neq \emptyset$ where $i\neq j$, $j\in [3]$. Without loss of generality, we assume $i$ = 1 and $j$ = 2. Let $x\in B_2$. Since $G$ is diamond-free, $x$ and $v_3$ are non-adjacent. If $N(x)\subseteq N(v_3)$, then the result follows from Lemma \ref{independent}. So assume there exists a vertex $y\in N(x)\setminus N(v_3)$. Then since $B_1 = \emptyset$, $y$ is not adjacent to $v_1$, so $\{v_1, v_3, x, y\}$ induces a 2$K_2$, a contradiction. This proves the Claim. \\

If there are adjacent vertices $x, y\in B_i$, for some $i\in$ [3], then $\{x, y, v_{i+1}, v_{i+2}\}$ induces a $2K_2$, a contradiction. Hence, $B_i$ is an independent set for all $i\in$ [3]. If $b_i \in B_i$ is adjacent to a vertex $z$ which is anticomplete to $Q$, then $\{b_i, z, v_{i+1}, v_{i+2}\}$ induces a $2K_2$, a contradiction. Therefore $V(G)$ = $Q \cup B_1 \cup B_2 \cup B_3$ and each $B_i$ is non-empty. We have two cases:\\

\noindent
\textit{Case 1}: \textbf{$|B_i|$ = 1 for some $i\in$ [3].} \\\noindent
Without loss of generality, let $B_1$ = $\{b_1\}$. Let $\alpha$ be any $\ell$-coloring of $G$. We prove that starting from $\alpha$, we can recolor the vertices of $G$ to a $\chi$-coloring, say $\alpha_1$, of $G$ such that $\alpha_1$ induces the color classes $\{v_1\}\cup B_3$, $\{v_2\}\cup B_1$, and $\{v_3\}\cup B_2$. Without loss of generality, let $\alpha(v_1)$ = 1, let $\alpha(v_2)$ = 2, and let $\alpha(v_3)$ = 3. Let $\alpha(b_1)$ = $r$, for some $r\in \{2,\dots,\ell\}$. 

Let $r\in \{2,3\}$. If $r$ = 2, i.e. $\alpha(b_1)$ = 2, then recolor every vertex of $B_2$ with color 3 and recolor every vertex of $B_3$ with color 1.

If $r$ = 3, i.e. $\alpha(b_1)$ = 3, then do the following:
\begin{itemize}
\setlength\itemsep{0em}
    \item recolor every vertex of $B_3$ with color 2,
    \item recolor every vertex of $B_2$ with color 4,
    \item recolor every vertex of $B_3$ with color 1,
    \item recolor $b_1$ with color 2,
    \item recolor every vertex of $B_2$ with color 3.
\end{itemize}
In each case, we obtain a $\chi$-coloring of $G$ isomorphic to $\alpha_1$ by recoloring each vertex at most twice. Then by the Renaming Lemma, we can reach the coloring $\alpha_1$ from $\alpha$ by recoloring each vertex at most twice. Therefore there exists a path between $\alpha$ and $\alpha_1$ in $R_{\ell}(G)$ where each vertex is recolored at most 4 $< n$ times. This implies $\alpha_1$ is a good coloring of $G$, a contradiction.

Let $r\notin \{2,3\}$. Without loss of generality, let $r$ = 4, i.e. $\alpha(b_1)$ = 4. Recolor every vertex of $B_2$ with color 3, recolor every vertex of $B_3$ with color 1, and recolor $b_1$ with color 2. We obtain a $\chi$-coloring of $G$ isomorphic to $\alpha_1$ by recoloring each vertex at once. Then by the Renaming Lemma, we can reach the coloring $\alpha_1$ from $\alpha$ by recoloring each vertex at most twice. Therefore, there exists a path between $\alpha$ and $\alpha_1$ where each vertex is recolored at most 3 $< n$ times. This implies $\alpha_1$ is a good coloring of $G$, a contradiction.\\

\noindent
\textit{Case 2}: \textbf{$|B_i| > 1$ for all $i$}.\\
\noindent
Without loss of generality, let $x_1$ and $x_2$ be in $B_1$. If $N(x_1)\subseteq N(x_2)$, then since $x_1$ and $x_2$ are non-adjacent, the result follows from Lemma \ref{independent}. So assume there exist two vertices $y_1\in N(x_1)$ and $z_2\in N(x_2)$, such that $y_1\notin N(x_2)$ and $z_2\notin N(x_1)$. If $y_1$ and $z_2$ are both in $B_2$ or $B_3$, then since the $B_i$s are an independent set, $\{x_1, y_1, x_2, z_2\}$ induces a $2K_2$, a contradiction. So without loss of generality, assume that $y_1\in B_2$ and $z_2\in B_3$. Since each $B_i$ has at least two vertices, there exist two vertices $y_2\in B_2$ and $z_1\in B_3$.

If $y_2$ is non-adjacent to both $x_2$ and $z_2$, then $\{x_2,z_2,y_2,v_2\}$ induces a $2K_2$, a contradiction. Without loss of generality, let $x_2y_2\in E(G)$. We make the following observations,
\begin{itemize}
\setlength\itemsep{0em}
    \item $y_1z_2\in E(G)$ since otherwise $\{x_1,y_1,x_2,z_2\}$ induces a $2K_2$.
    \item $x_1y_2\in E(G)$ since otherwise $\{x_1,y_1,x_2,y_2\}$ induces a $2K_2$.
    \item $y_2z_2\in E(G)$ since otherwise $\{x_1,y_2,z_2,v_3\}$ induces a $2K_2$.
    \item $z_1$ is adjacent to $x_1$ or $y_1$ since otherwise $\{x_1,y_1,z_1,v_3\}$ induces a $2K_2$.
    \item If $z_1$ is adjacent to $x_1$, then
    \begin{itemize}
    \setlength\itemsep{0em}
          \item $z_1x_2\in E(G)$ since otherwise $\{x_1,z_1,x_2,z_2\}$ induces a $2K_2$
          \item $z_1y_1\in E(G)$ since otherwise $\{x_2,z_1,y_1,v_2\}$ induces a $2K_2$
          \item $z_1y_2\notin E(G)$ since otherwise $\{x_1,x_2,y_2,z_1\}$ induces a diamond.
    \end{itemize}
    \item If $z_1$ is adjacent to $y_1$, then
    \begin{itemize}   
    \setlength\itemsep{0em}
          \item {$z_1x_2\in E(G)$ since otherwise $\{y_1,z_1,x_2,v_1\}$ induces a $2K_2$}
          \item {$z_1y_2\notin E(G)$ since otherwise $\{x_2,z_2,y_2,z_1\}$ induces a diamond}
          \item {$z_1x_1\in E(G)$ since otherwise $\{x_1,y_2,z_1,v_3\}$ induces a $2K_2$.}
    \end{itemize}
\end{itemize}
Therefore $\{v_1,v_2,v_3,x_1,x_2,y_1,y_2,z_1,z_2\}$ induces a 3-prism star. We claim that $G$ is isomorphic to the 3-prism star. If not, then since $V(G)$ = $Q \cup B_1 \cup B_2 \cup B_3$ there exists a vertex $u\in B_i$ for some $i\in [3]$. Without loss of generality, let $u\in B_1$. The vertex $u$ is either complete or anticomplete to $\{y_2, z_1\}$; since otherwise $\{u, y_2, z_1, v_2\}$ or $\{u, y_2, z_1, v_3\}$ induces a $2K_2$, a contradiction. First assume $u$ is anticomplete to $\{y_2, z_1\}$. If $uz_2\notin E(G)$, then $\{y_2, z_2, u, v_1\}$ induces a $2K_2$, a contradiction. If $uz_2\in E(G)$, then $\{u, z_2, x_1, z_1\}$ induces a $2K_2$, a contradiction. Now assume $u$ is complete to $\{y_2, z_1\}$. If $uy_1\in E(G)$ ($uz_2 \in E(G)$), then $\{u, y_1, z_1, x_1\}$ ($\{u, y_2, z_2, x_2\}$) induces a diamond, a contradiction. But $uy_1\notin E(G)$ and $uz_2\notin E(G)$ implies $\{y_1, z_2, u, v_1\}$ induces a $2K_2$, a contradiction. Therefore, $G$ is isomorphic to the 3-prism star. From Lemma \ref{prism}, $G$ is good, a contradiction.
\end{proof}

\begin{theorem}\label{2k2diamond}
Every $(2K_2, diamond)$-free graph $G$ is recolorable with $\ell$ recoloring diameter at most 2$n^2$, for all $\ell\geq \chi(G)$+1.
\end{theorem}
\begin{proof}
We prove that every (2$K_2$, diamond)-free graph $G$ is good. If $G$ is triangle-free, then the result follows from Theorem \ref{paw2k2}. We assume that $G$ is connected, contains a triangle, and is not complete. The proof is by induction on $n$.

Let $Q$ = $\{v_1, v_2,\dots, v_p\}$, $p\geq 3$, be a maximal clique in $G$ with at least 3 vertices. For all $i\in [p]$, define $B_i$ = $\{v\in V(G)\mid N(v) \cap Q$ = $\{v_i\}\}$. Since $G$ is diamond-free, every neighbor of $v_i$ not in $Q$ is in $B_i$. Since $G$ is connected and not complete, there is some $B_i \neq \emptyset$ where $i\in [p]$. Without loss of generality, let $b_1\in B_1$. If $N(b_1)$ = $\{v_1\}$, then $N(b_1)\subseteq N(v_2)$. By the induction hypothesis, $G$-$b_1$ is good, and the result follows from Lemma \ref{independent}. So we assume $b_1$ is adjacent to a vertex $x\neq v_1$. If $p\geq 4$, then there exist unique vertices $v_j$ and $v_k$ such that $x \notin N(v_j)\cup N(v_k)$, where $j$, $k \in [p]$. This implies $\{v_j, v_k, b_1, x\}$ induces a $2K_2$, a contradiction. Therefore $p$ = 3 and thus the clique number of $G$ is 3. The result follows from Lemma \ref{clique3}.
\end{proof}

\begin{figure}[H]
\centering
\begin{tikzpicture}[scale=3]
\tikzstyle{vertex}=[circle, draw, fill=black, inner sep=0pt, minimum size=5pt]
    
    \draw (0,0) circle (1);
    
    \node[vertex, label=above:1](0) at (0,1) {};
    \node[vertex, label=above:2](1) at (cos{67.5},sin{67.5}) {};
    \node[vertex, label=above:3](2) at (cos{45},sin{45}) {};
    \node[vertex, label=above:4](3) at (cos{22.5}, sin{22.5}) {};
    \node[vertex, label=right:5](4) at (1,0) {};
    \node[vertex, label=below:7](5) at (cos{22.5}, -sin{22.5}) {};
    \node[vertex, label=below:2](6) at (cos{45},-sin{45}) {};
    \node[vertex, label=below:3](7) at (cos{67.5},-sin{67.5}) {};
    \node[vertex, label=below:4](8) at (0,-1) {};
    \node[vertex, label=left:5](9) at (-cos{67.5},-sin{67.5}) {};
     \node[vertex, label=left:1](10) at (-cos{45},-sin{45}) {};
    \node[vertex, label=left:2](11) at (-cos{22.5}, -sin{22.5}) {};
    \node[vertex, label=left:3](12) at (-1,0) {};
    \node[vertex, label=above:4](13) at (-cos{22.5},sin{22.5}) {};
    \node[vertex, label=above:6](14) at (-cos{45},sin{45}) {};
    \node[vertex, label=above:7](15) at (-cos{67.5},sin{67.5}) {};
    
    \draw(0)--(2);
    \draw(0)--(3);
    \draw(0)--(4);
    \draw(0)--(5);
    \draw(0)--(6);
    \draw(0)--(7);
    \draw(0)--(9);
    \draw(0)--(11);
    \draw(0)--(12);
    \draw(0)--(13);
    \draw(0)--(14);
    \draw(1)--(3);
    \draw(1)--(4);
    \draw(1)--(5);
    \draw(1)--(7);
    \draw(1)--(8);
    \draw(1)--(10);
    \draw(1)--(12);
    \draw(1)--(13);
    \draw(1)--(14);
    \draw(1)--(15);
    \draw(2)--(4);
    \draw(2)--(5);
    \draw(2)--(6);
    \draw(2)--(8);
    \draw(2)--(9);
    \draw(2)--(11);
    \draw(2)--(13);
    \draw(2)--(14);
    \draw(2)--(15);
    \draw(3)--(5);
    \draw(3)--(6);
    \draw(3)--(7);
    \draw(3)--(9);
    \draw(3)--(10);
    \draw(3)--(12);
    \draw(3)--(14);
    \draw(3)--(15);
    \draw(4)--(6);
    \draw(4)--(7);
    \draw(4)--(8);
    \draw(4)--(10);
    \draw(4)--(11);
    \draw(4)--(13);
    \draw(4)--(14);
    \draw(4)--(15);
    \draw(5)--(7);
    \draw(5)--(8);
    \draw(5)--(9);
    \draw(5)--(10);
    \draw(5)--(11);
    \draw(5)--(12);
    \draw(5)--(14);
    \draw(6)--(8);
    \draw(6)--(9);
    \draw(6)--(10);
    \draw(6)--(12);
    \draw(6)--(13);
    \draw(6)--(15);
    \draw(7)--(9);
    \draw(7)--(10);
    \draw(7)--(11);
    \draw(7)--(13);
    \draw(7)--(14);
    \draw(8)--(10);
    \draw(8)--(11);
    \draw(8)--(12);
    \draw(8)--(14);
    \draw(8)--(15);
    \draw(9)--(11);
    \draw(9)--(12);
    \draw(9)--(13);
    \draw(9)--(14);
    \draw(9)--(15);
    \draw(10)--(12);
    \draw(10)--(13);
    \draw(10)--(14);
    \draw(10)--(15);
    \draw(11)--(13);
    \draw(11)--(14);
    \draw(11)--(15);
    \draw(12)--(14);
    \draw(12)--(15);
    \draw(13)--(15);

\end{tikzpicture}
\hspace{0mm}
\begin{tikzpicture}[scale=3]
\tikzstyle{vertex}=[circle, draw, fill=black, inner sep=0pt, minimum size=5pt]
    
    \draw (0,0) circle (1);
    
    \node[vertex, label=above:1](0) at (0,1) {};
    \node[vertex, label=above:2](1) at (cos{67.5},sin{67.5}) {};
    \node[vertex, label=above:3](2) at (cos{45},sin{45}) {};
    \node[vertex, label=above:4](3) at (cos{22.5}, sin{22.5}) {};
    \node[vertex, label=right:5](4) at (1,0) {};
    \node[vertex, label=below:6](5) at (cos{22.5}, -sin{22.5}) {};
    \node[vertex, label=below:7](6) at (cos{45},-sin{45}) {};
    \node[vertex, label=below:8](7) at (cos{67.5},-sin{67.5}) {};
    \node[vertex, label=below:1](8) at (0,-1) {};
    \node[vertex, label=left:2](9) at (-cos{67.5},-sin{67.5}) {};
     \node[vertex, label=left:3](10) at (-cos{45},-sin{45}) {};
    \node[vertex, label=left:4](11) at (-cos{22.5}, -sin{22.5}) {};
    \node[vertex, label=left:5](12) at (-1,0) {};
    \node[vertex, label=above:6](13) at (-cos{22.5},sin{22.5}) {};
    \node[vertex, label=above:7](14) at (-cos{45},sin{45}) {};
    \node[vertex, label=above:8](15) at (-cos{67.5},sin{67.5}) {};
    
    \draw(0)--(2);
    \draw(0)--(3);
    \draw(0)--(4);
    \draw(0)--(5);
    \draw(0)--(6);
    \draw(0)--(7);
    \draw(0)--(9);
    \draw(0)--(11);
    \draw(0)--(12);
    \draw(0)--(13);
    \draw(0)--(14);
    \draw(1)--(3);
    \draw(1)--(4);
    \draw(1)--(5);
    \draw(1)--(7);
    \draw(1)--(8);
    \draw(1)--(10);
    \draw(1)--(12);
    \draw(1)--(13);
    \draw(1)--(14);
    \draw(1)--(15);
    \draw(2)--(4);
    \draw(2)--(5);
    \draw(2)--(6);
    \draw(2)--(8);
    \draw(2)--(9);
    \draw(2)--(11);
    \draw(2)--(13);
    \draw(2)--(14);
    \draw(2)--(15);
    \draw(3)--(5);
    \draw(3)--(6);
    \draw(3)--(7);
    \draw(3)--(9);
    \draw(3)--(10);
    \draw(3)--(12);
    \draw(3)--(14);
    \draw(3)--(15);
    \draw(4)--(6);
    \draw(4)--(7);
    \draw(4)--(8);
    \draw(4)--(10);
    \draw(4)--(11);
    \draw(4)--(13);
    \draw(4)--(14);
    \draw(4)--(15);
    \draw(5)--(7);
    \draw(5)--(8);
    \draw(5)--(9);
    \draw(5)--(10);
    \draw(5)--(11);
    \draw(5)--(12);
    \draw(5)--(14);
    \draw(6)--(8);
    \draw(6)--(9);
    \draw(6)--(10);
    \draw(6)--(12);
    \draw(6)--(13);
    \draw(6)--(15);
    \draw(7)--(9);
    \draw(7)--(10);
    \draw(7)--(11);
    \draw(7)--(13);
    \draw(7)--(14);
    \draw(8)--(10);
    \draw(8)--(11);
    \draw(8)--(12);
    \draw(8)--(14);
    \draw(8)--(15);
    \draw(9)--(11);
    \draw(9)--(12);
    \draw(9)--(13);
    \draw(9)--(14);
    \draw(9)--(15);
    \draw(10)--(12);
    \draw(10)--(13);
    \draw(10)--(14);
    \draw(10)--(15);
    \draw(11)--(13);
    \draw(11)--(14);
    \draw(11)--(15);
    \draw(12)--(14);
    \draw(12)--(15);
    \draw(13)--(15);
    
\end{tikzpicture}
\caption{A 7-coloring and a frozen 8-coloring of $G$ \cite{2K2free}.}
\label{fig:G}
\end{figure}
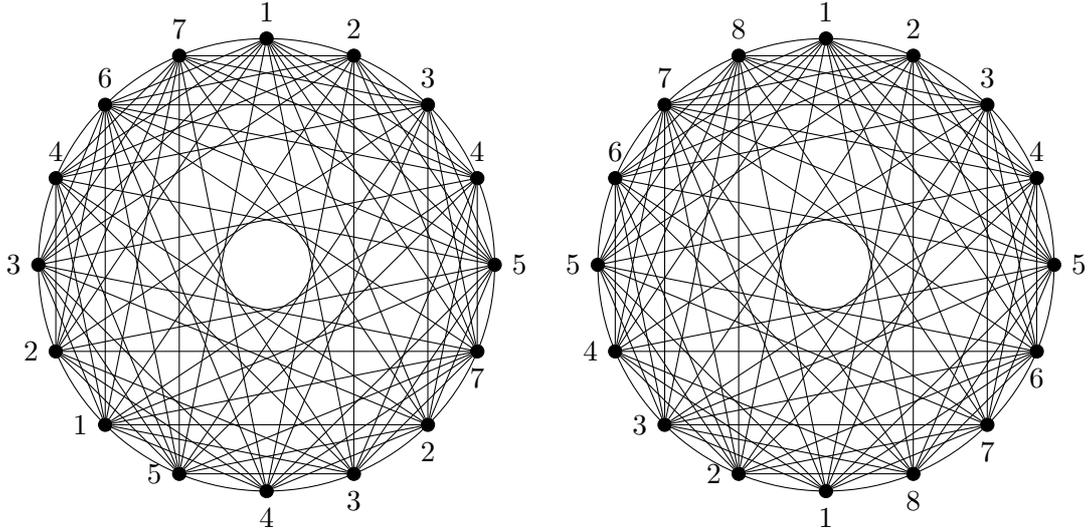

\begin{theorem}
For all $p\geq 1$, there exists a $k$-colorable $(2K_2, 4K_1, co$-$diamond, co$-$claw)$-free graph that is not $(k+ p)$-mixing.
\end{theorem}
\begin{proof}
Let $G$ be the graph of Figure \ref{fig:G}. Then $G$ is a 7-colorable (2$K_2$, 4$K_1$, co-diamond, co-claw)-free graph which admits a frozen 8-coloring \cite{2K2free}. For $p\geq 2$, we take the pairwise join of $p$ copies of $G$ and the resulting graph will be a $7p$-colorable (2$K_2$, 4$K_1$, co-diamond, co-claw)-free graph which admits a frozen $8p$-coloring.
\end{proof}


Jamison and Olariu gave a forbidden induced subgraph characterization for the class of $P_4$-sparse graphs. 

\begin{theorem}(\cite{Jamison})
    A graph $G$ is $P_4$-sparse if and only if $G$ is ($P_5$, $C_5$, house, fork, co-fork, banner, co-banner)-free.
\end{theorem}

\begin{figure}[h]
    \centering
    \begin{tikzpicture}[scale=0.8]
    \tikzstyle{vertex}=[circle, draw, fill=black, inner sep=0pt, minimum size=3.75pt]
    \node[vertex](1) at (-2,0.15){};
    \node[vertex](2) at (-1,0.15){};
    \node[vertex](3) at (-2,1){};
    \node[vertex](4) at (-1,1){};
    \node[vertex](5) at (-1.5,1.75){};
    \draw(1)--(2);
    \draw(4)--(2);
    \draw(3)--(4);
    \draw(1)--(3);
    \draw(4)--(5);
    \draw(3)--(5);
    \node[] at (-1.5, -1){$house$};    
    \node[vertex](6) at (0.5,1){};
    \node[vertex](7) at (1.5,1){};
    \node[vertex](8) at (2.5,1){};
    \node[vertex](9) at (3.5,1.5){};
    \node[vertex](10) at (3.5,0.5){};
    \node[] at (1.7, -1){$fork$};  
    \draw(6)--(7);
    \draw(8)--(7);
    \draw(8)--(9);
    \draw(8)--(10);
    \node[vertex](11) at (5,1){};
    \node[vertex](12) at (6,1.5){};
    \node[vertex](13) at (6,0.5){};
    \node[vertex](14) at (7,1){};    
    \node[vertex](15) at (8,1){};
    \draw(12)--(11);
    \draw(13)--(11);
    \draw(12)--(13);
    \draw(12)--(14);
    \draw(14)--(13);
    \draw(14)--(15);
    \node[] at (6.2, -1){$co$-$fork$};
    \node[vertex](16) at (9.5,1.5){};
    \node[vertex](17) at (10.5,1.5){};
    \node[vertex](18) at (9.5,0.5){};
    \node[vertex](19) at (10.5,0.5){};
    \node[vertex](20) at (11.5,0.5){};
    \node[] at (10.25, -1){$banner$};
    \draw(18)--(16);
    \draw(16)--(17);
    \draw(17)--(19);
    \draw(18)--(19);
    \draw(19)--(20);
    \node[vertex](21) at (13,1.5){};
    \node[vertex](22) at (13,0.5){};
    \node[vertex](23) at (14,1){};
    \node[vertex](24) at (15,1){};
    \node[vertex](25) at (16,1){};
    \node[] at (14.5, -1){$co$-$banner$};
    \draw(21)--(22);
    \draw(21)--(23);
    \draw(22)--(23);
    \draw(23)--(24);
    \draw(24)--(25);
    \end{tikzpicture}
\end{figure}

Therefore the class of ($P_5$, $C_5$,\ house,\ co-banner)-free graphs strictly contains the class of $P_4$-sparse graphs. In \cite{feghali2020}, Feghali and Fiala proved that every 3-colorable ($P_5$, $C_5$, house)-free graph is recolorable. They also asked if the diameter of $R_{k+1}(G)$ is finite for a $k$-colorable ($P_5$, $C_5$, house)-free graph $G$.

\begin{theorem}\label{P4sparselemma}
    Every ($P_5$, $C_5$,\ house,\ co-banner)-free graph $G$ is recolorable with $\ell$ recoloring diameter at most 2$n^2$, for all $\ell\geq \chi(G)$+1.
\end{theorem}
\begin{proof}
We prove that every ($P_5$, $C_5$,\ house, co-banner)-free graph $G$ is good. The proof is by induction on $n$. If $G$ is $P_4$-free, then either $G$ is a disjoint union or a join of two graphs and the result follows from Lemma \ref{join}. So we assume that $G$ contains an induced $P_4$. Let $\{u, v, x, y\} \subseteq V(G)$ induce a $P_4$, where $\{uv, vx, xy\}\subseteq E(G)$.

\textbf{Claim}: $N(y)\subseteq N(v)$.\\
If not, then there exists a vertex $z \in N(y)\setminus N(v)$. If $z$ is anticomplete to $\{x,u\}$, then $\{u, v, x, y, z\}$ induces a $P_5$, a contradiction. If $z$ is adjacent to $x$ but not to $u$, then $\{u, v, x, y, z\}$ induces a co-banner, a contradiction. If $z$ is adjacent to $u$ but not to $x$, then $\{u, v, x, y, z\}$ induces a $C_5$, a contradiction. If $z$ is complete to $\{x,u\}$, then $\{u, v, x, y, z\}$ induces a house, a contradiction. This proves the Claim.

Therefore $G$ contains two independent vertices $y$ and $v$, such that $N(y)\subseteq N(v)$. By the induction hypothesis, $G$-$y$ is good and the result follows from Lemma \ref{independent}.
\end{proof}

\section{Conclusion}

In \cite{Belavadi}, it was proved that every $H$-free graph $G$ is recolorable if and only if $H$ is an induced subgraph of $P_4$ or $P_3$+$P_1$. Let $H_1$ and $H_2$ be two 4-vertex graphs different from $P_4$ and $P_3$+$P_1$. As explained in the introduction, if every ($H_1$, $H_2$)-free graph $G$ is recolorable, then either $H_1$ or $H_2$ is isomorphic to 2$K_2$. We presented several recolorability results for the class of (2$K_2$, $H$)-free graphs, where $H$ is a 4-vertex graph. To obtain a dichotomy for (2$K_2$, $H$)-free graphs related to recolorability, one needs to know whether every (2$K_2$, $K_4$)-free graph is recolorable. We proved that every (2$K_2$, triangle)-free graph is recolorable. It remains open to determine whether every (2$K_2$, $K_4$)-free graph containing a triangle is recolorable.

\end{document}